\documentclass[12pt,twoside]{amsart}
\usepackage{amssymb}
\usepackage{graphicx}
\usepackage{amsmath}

\newtheorem{theorem}{Theorem}[section]
\newtheorem{corollary}[theorem]{Corollary}

\newtheorem{example}[theorem]{Example}

\newtheorem{lemma}[theorem]{Lemma}

\newtheorem{proposition}[theorem]{Proposition}

\begin{document}
\begin{center}
\textbf{{\Large New constructions of two slim dense near
hexagons}}

\vskip .6cm

{\large Binod Kumar Sahoo}\smallskip\\
{\small Department of Mathematics\\
National Institute of Technology\\
Rourkela-769008, Orissa, India.\\
E-mail: binodkumar@gmail.com}\\
(05 August 2006)
\end{center}

\begin{abstract}
We provide a geometrical construction of the slim dense near hexagon
with parameters $(s,t,t_{2})=(2,5,\{1,2\})$. Using this
construction, we construct the rank 3 symplectic dual polar space
$DSp(6,2)$ which is the slim dense near hexagon with parameters
$(s,t,t_{2})=(2,6,2)$. Both the near hexagons are constructed from
two copies of a generalized quadrangle with parameters (2,2).\\

\end{abstract}

\title[Constructions of slim dense near hexagons]{}
\author[B. K. Sahoo]{}
\address[]{}
\maketitle

\section{Introduction}

A \textit{partial linear space} is a point-line geometry $S=(P,L)$
with `point-set' $P$ and `line-set' $L$ of subsets of $P$ of size
at least 2 such that any two distinct points are contained in at
most one line. Two distinct points $x$ and $y$ are
\textit{collinear}, written as $x\sim y$, if there is a line
containing them. In that case we denote this line by $xy$. We
write $x\nsim y$ if two distinct points $x$ and $y$ are not
collinear. For $ x\in P$ and $A\subseteq P$, we define $x^{\perp
}=\{x\}\cup \{y\in P:x\sim y\}$ and $A^{\perp }=\underset{x\in
A}{\cap }x^{\perp }$. Note that $xy=\{x,y\}^{\perp}$ if $x\sim y$.
A subset of $P$ is a \textit{subspace} of $S$ if any line
containing at least two of its points is contained in it. A
\textit{geometric hyperplane} of $S$ is a subspace, different from
the empty set and $P$, that meets every line non-trivially. The
graph $\Gamma (P)$ with vertex set $P$, in which two distinct
vertices are {\it adjacent} if they are collinear in $S$, is the
\textit{collinearity graph} of $S$. The diameter of $\Gamma (P)$
is called the \textit{diameter} of $S$. If $\Gamma (P)$ is
connected, we say that $S$ is {\it connected}.

A \textit{near polygon} \cite{SY} is a connected partial linear
space of finite diameter satisfying the following near-polygon
property $(NP)$: \textit{For each point-line pair $(x,l)$ with
$x\notin l$, there is a unique point in $l$ nearest to $x$.} If the
diameter of $S$ is $n$, then $S$ is called a \textit{near $2n$-gon},
and for $n=3$ a \textit{near hexagon}. A near 0-gon is a point and a
near 2-gon is a line. An important class of near polygons is the the
class of generalized $2n$-gons. Near 4-gons are precisely the class
of \textit{generalized quadrangles}.

Let $S=(P,L)$ be a near polygon. A subset $C$ of $P$ is
\textit{convex} if every shortest path in $\Gamma (P)$ between two
points of $C$ is entirely contained in $C$. A \textit{quad} is a
convex subset of $P$ of diameter $2$ such that no point of it is
adjacent to all other points of it. We denote by $d(x,y)$ the
\textit{distance} in $\Gamma (P)$ between two points $x$ and $y$ of
$S$. Let $x_{1}$ and $x_{2}$ be two points of $S$ with
$d(x_{1},x_{2})=2$. If $x_{1}$ and $x_{2}$ have at least two common
neighbours $y_{1}$ and $y_{2}$ such that one of the lines
$x_{i}y_{j}$ contains at least three points, then $x_{1}$ and
$x_{2}$ are contained in a unique quad which is a generalized
quadrangle (\cite{SY}, Proposition 2.5, p.10). A near polygon is
called \textit{dense} if every pair of points at distance 2 have at
least two common neighbours. A structure theory for dense near
polygons can be found in \cite{BW}. Dual polar spaces were
characterized in terms of point-line geometries by Cameron. We refer
to \cite{Ca} for the definition and the related results concerning
dual polar space. The class of dual polar spaces coincides with the
class of classical dense near polygons.

In a dense near polygon, the number of lines through a point is
independent of that point (\cite{BW}, Lemma 19, p.152). A dense near
$2n$-gon is said to have parameters $(s,t)$ if each line contains
$s+1$ points and each point is contained in $t+1$ lines. For $n=2$,
a generalized quadrangle with parameters $(s,t)$ is written as
$(s,t)$-GQ, for short. A near polygon is \textit{slim} if each line
contains exactly three points. In that case if $x$ and $y$ are two
collinear points, we define $x*y$ by $xy=\{x,y,x*y\}$. All slim
dense near hexagons have been classified by Brouwer et al.
\cite{BCHW}. We refer to \cite{Bart-book} and \cite{V} for some new
classification results concerning slim dense near polygons.

\begin{theorem}\label{classification-result}
(\cite{BCHW}, Theorem 1.1, p.349) Any slim dense near hexagon is
finite and is isomorphic to one of the eleven near hexagons with
parameters as given below:
\begin{center}
$
\begin{array}{|c|c|c|c|c|c|c|c|c|c|c|c|}\hline
v & 891 & 759 & 729 & 567 & 405 & 243 & 135 & 105 & 81 & 45 & 27
\\\hline
t & 20 & 14 & 11 & 14 & 11 & 8 & 6 & 5 & 5 & 3 & 2 \\\hline t_{2} &
4 & 2 & 1 & 2,4 & 1,2,4 & 1,4 & 2 & 1,2 & 1,4 & 1,2 & 1\\ \hline
\end{array} $
\end{center}
\end{theorem}

Here $v$ is the number of points, $t$ is the number of lines
containing a point and $t_{2}$ is defined by
$t_{2}+1=|x^{\perp}\cap y^{\perp}|$ for points $x,y$ with
$d(x,y)=2$ (though this depends on $x$ and $y$).

The slim dense near hexagon on 135 points with parameters
$(s,t,t_{2})=(2,6,2)$ is the symplectic dual polar space of rank 3
over the field $F_{2}$ with two elements. We denote this near
hexagon by $DSp(6,2)$ and by $\mathbb{H}_{3}$ the near hexagon on
105 points with parameters $(s,t,t_{2})=(2,5,\{1,2\})$. In the
following, we present three known descriptions of the near hexagon
$\mathbb{H}_{3}$. An equivalent description of this near hexagon is
given in (\cite{BCHW}, p.352) using the notion of a Fischer space.

\begin{example}
(\cite{BCHW}, p.355) Let $X$ be a set of size 8. Let $S$ be the
partial linear space whose point set consists of all partitions of
$X$ into four 2-subsets, and lines are the collections of partitions
sharing two 2-subsets. Then $S$ is a slim dense near hexagon
isomorphic to $\mathbb{H}_{3}$.
\end{example}

\begin{example}\label{example2}
(\cite{BCHW}, p.352) Let $Q(6,2)$ be a non-singular quadric in
$PG(6,2)$ and let $DQ(6,2)$ be the associated dual polar space.
Let $\Pi$ be a hyperplane of $PG(6,2)$ intersecting $Q(6,2)$ in a
non-singular hyperbolic quadric $Q^{+}(5,2)$. The set of all
maximal subspaces of $Q(6,2)$ which are not contained in
$Q^{+}(5,2)$ is a subspace of $DQ(6,2)$ isomorphic to the near
hexagon $\mathbb{H}_{3}$. (Note that $DQ(6,2)$ is isomorphic to
$DSp(6,2)$.)
\end{example}

The next construction (Example \ref{bart-de-bruyn}) is due to Bart
De Bruyn (\cite{D}, p.51). We first note down some facts about
$(2,2)$-GQs which play an important role in Example
\ref{bart-de-bruyn} below as well as in our constructions in the
next section.

Let $S=(P,L)$ be a $(2,2)$-GQ. Then $S$ is unique up to an
isomorphism (\cite{PT}, 5.2.3, p.78). In the notation of
\cite{PT}, $S$ is isomorphic to the classical generalized
quadrangle $W(2)\simeq Q(4,2)$. \textbf{A combinatorial model of a
$(2,2)$-GQ:} Let $\Omega =\{1,2,3,4,5,6\}$. An \textit{edge} of
$\Omega $ is a $2$-subset of $ \Omega $. A \textit{factor}\ of
$\Omega $ is a set of three pair-wise disjoint edges. Let
$\mathcal{E}$ be the set of all edges and $\mathcal{F}$ be the set
of all factors of $\Omega$. Then $|\mathcal{E}|=|\mathcal{F}|=15$.
Taking $\mathcal{E}$ as the point-set and $\mathcal{F}$ as the
line-set, the pair $(\mathcal{E},\mathcal{F})$ is a $(2,2)$-GQ.

A $k$\textit{-arc} of points of $S$ is a set of $k$ pair-wise
non-collinear points of $S$. A $k$-arc is \textit{complete} if it is
not contained in a $(k+1)$-arc. A point $x$ is a \textit{center} of
a $k$-arc if $x$ is collinear with every point of it. A 3-arc is
called a \textit{triad}. We refer to \cite{S} for a detailed study
of $k$-arcs of a $(2,t)$-GQ. We repeatedly use the following result,
mostly without mention.

\begin{lemma}\label{facts-triads}
(\cite{S}, Proposition 3.2, p.161) Let $S=(P,L)$ be a (2,2)-GQ and
$T$ be a triad of points of $S$. Then
\begin{enumerate}
\item[$(i)$] $|T^{\perp}|=1$ or $3$;

\item[$(ii)$] $|T^{\perp }|=1$ if and only if $T$ is an incomplete
triad if and only if $T$ is contained in a unique (2,1)-subGQ of
$S$;

\item[$(iii)$] $|T^{\perp}| =3$ if and only if $T$ is a complete
triad.
\end{enumerate}
\end{lemma}

Dually, we can define a $k$-arc of lines of $S$ and center of such
a $k$-arc. Since $W(2)$ is self-dual (\cite{PT}, 3.2.1, p.43),
Lemma \ref{facts-triads} holds for a triad of lines also. From the
combinatorial model above, it can be seen that any two
non-collinear points (respectively; disjoint lines) of $S$ are
contained in a unique complete triad of points (respectively;
lines) of $S$.

\begin{example}\label{bart-de-bruyn}
(\cite{D}, p.51) Let $S=(P,L)$ be the $(2,2)$-GQ. A partial linear
space $S'=(P',L')$ can be constructed from $S$ as follows. Set
$$P'=\{(x,y)\in P\times P:x=y \text{
or }x\sim y\}.$$ There are four types of lines in $L'$. For each
line $l=\{x,y,z\}$ of $S$, we have the following lines of $S'$ (up
to a permutation of the points of $l$):
\begin{enumerate}
\item[$(i)$] $\{(x,x),(y,y),(z,z)\}$,

\item[$(ii)$] $\{(x,x),(x,y),(x,z)\}$,

\item[$(iii)$] $\{(x,y),(y,z),(z,x)\}$.
\end{enumerate}
Let $T=\{k,m,n\}$ be an incomplete triad of lines of $S$ such that
$T^{\perp}=\{l\}$. We may assume that $k\cap l=\{x\}$, $m\cap
l=\{y\}$ and $n\cap l=\{z\}$. Let $x'$ be an arbitrary point of
$k\setminus \{x\}$. Let $y'$ (respectively, $z'$) be the unique
point of $m\setminus \{y\}$ (respectively, $n\setminus \{z\}$) not
collinear with $x'$. Then the lines of fourth type are:
\begin{enumerate}
\item[$(iv)$] $\{(x,x'),(y,y'),(z,z')\}$.
\end{enumerate}
Then $S'$ is a slim dense near hexagon isomorphic to
$\mathbb{H}_{3}$. (It can be seen that the set $\{x',y',z'\}$ is an
incomplete triad of points of $S$.)
\end{example}

\section{New Constructions}

Here, we provide new geometrical constructions for the near
hexagons $DSp(6,2)$ and $\mathbb{H}_{3}$ from two copies of a
(2,2)-GQ. In fact, we first construct $\mathbb{H}_{3}$ and then
construct $DSp(6,2)$ in which $\mathbb{H}_{3}$ is embedded as a
geometric hyperplane.

\subsection{Construction of $\mathbb{H}_{3}$}\label{H_3}

In this construction the point set is the same as in Example
\ref{bart-de-bruyn} (up to an isomorphism). Let $S=(P,L)$ and
$S'=(P',L')$ be two (2,2)-GQs. Let $x\leftrightarrow x'$, $x\in
P,x'\in P'$, be an isomorphism between $S$ and $S'$. We define a
partial linear space $\mathcal{S}=(\mathcal{P},\mathcal{L})$ as
follows. The point set is
$$\mathcal{P}=\{(x,y')\in P\times P': y'\in x'^{\perp}\},$$
and the lines are of the form $$\{(x,u'),(y,v'),(z,w')\},$$ where
$T=\{x,y,z\}$ is either a line or a complete triad and
$T'^{\perp}=\{u',v',w'\}$.

\begin{theorem}\label{first-main-theorem}
The partial linear space $\mathcal{S}=(\mathcal{P},\mathcal{L})$
is a slim dense near hexagon with parameters
$(s,t,t_{2})=(2,5,\{1,2\})$.
\end{theorem}

\subsection{Construction of $DSp(6,2)$}

Let $S=(P,L)$, $S'=(P',L')$ and
$\mathcal{S}=(\mathcal{P},\mathcal{L})$ be as in the construction of
$\mathbb{H}_{3}$. We define a partial linear space
$\mathbb{S}=(\mathbb{P},\mathbb{L})$ as follows. The point set is
$$\mathbb{P}=\mathcal{P}\cup P\cup P',$$
and the line set $\mathbb{L}$ consists of the lines in $\mathcal{L}$
together with the following collection $\mathbb{L}_{1}$ of lines:
\begin{enumerate}
\item[$(\mathbb{L}_{1}):$] The collection $\mathbb{L}_{1}$
consists of lines of the form $\{x,(x,u'),u'\}$ for every point
$(x,u')\in \mathcal{P}$.
\end{enumerate}

\begin{theorem}\label{second-main-theorem}
The partial linear space $\mathbb{S}=(\mathbb{P},\mathbb{L})$ is a
slim dense near hexagon with parameters $(s,t,t_{2})=(2,6,2)$
\end{theorem}

An immediate consequence of Theorems \ref{first-main-theorem} and
\ref{second-main-theorem} is that the near hexagon $\mathbb{H}_{3}$
is a geometric hyperplane of the near hexagon $DSp(6,2)$.

\section{Proof of Theorem \ref{first-main-theorem}}
\label{Proof-of-first-main-theorem}

Let $\alpha = (x,u')$ and $\beta = (y,v')$ be two distinct points of
$\mathcal{S}$. By the construction of lines of $\mathcal{S}$,
$\alpha\sim \beta$ if and only if $x\neq y$, $u'\neq v'$, $u'\in
y'^{\perp}$ and $v'\in x'^{\perp}$. Let $\alpha $ and $\beta $ be
distinct non-collinear points of $\mathcal{S}$. Then one of the
following possibilities occur:
\begin{enumerate}
\item[$(A1)$] $x=y,u'\neq v'$;

\item[$(A2)$] $x\neq y,u'=v'$;

\item[$(A3)$] $x\neq y,u'\neq v'$, $u'\notin y'^{\perp}$ and $v'\notin x'^{\perp}$;

\item[$(A4)$] $x\neq y,u'\neq v'$ and either $u'\in y'^{\perp}$ and $v'\notin x'^{\perp}$
or $u'\notin y'^{\perp}$ and $v'\in x'^{\perp}$.
\end{enumerate}

\begin{lemma}\label{one-argument-coincide}
Assume that $(A1)$ or $(A2)$ holds. Then
$|\{\alpha,\beta\}^{\perp}|\geq 2$.
\end{lemma}

\begin{proof}
Assume that $(A1)$ holds. Then $x'\in \{u',v'\}^{\perp}$. If $u'\sim
v'$, we may assume that $v'\neq x'$. So $x'v'=u'v'$. Then
$(v,u'*v')$ and $(x*v,u'*v')$ are in $\{\alpha,\beta\}^{\perp}$. If
$u'\nsim v'$, let $\{u',v',w'\}$ be the complete triad of $S'$
containing $u'$ and $v'$. Then $(a,w')$ and $(b,w')$ are in
$\{\alpha,\beta\}^{\perp}$, where $\{a',b'\}=
\{u',v'\}^{\perp}\setminus \{x'\}$.

Now assume that $(A2)$ holds. Then $u\in \{x,y\}^{\perp}$ in $S$. If
$x\sim y$, we assume that $u\neq y$. Then $(x*y,y')$ and $(x*y,
y'*u')$ are in $\{\alpha,\beta\}^{\perp}$. If $x\nsim y$, let
$\{x,y,w\}$ be the complete triad of $S$ containing $x$ and $y$. Let
$\{x,y\}^{\perp}\setminus \{u\}=\{a,b\}$ in $S$. Then $(w,a')$ and
$(w,b')$ are in $\{\alpha,\beta\}^{\perp}$.
\end{proof}

\begin{lemma}\label{no-argument-coincide-distance-2}
Assume that $(A3)$ holds. Then $|\{\alpha,\beta\}^{\perp}|\geq 3$.
\end{lemma}

\begin{proof}
If $x\sim y$ and $u'\sim v'$, then $\{x',y',u',v'\}$ form a
quadrangle in $S'$. Then $(u,x')$, $(v,y')$ and $(u*v,x'*y')$ are
in $\{\alpha,\beta\}^{\perp}$.

If $x\sim y$ and $u'\nsim v'$, let $T'=\{u'*x',v'*y',z'\}$ be the
complete triad of $S'$ containing $u'*x'$ and $v'*y'$. Then
$u',v'\in T'^{\perp}$ and $x'*y'\notin T'$. Now, $x'*y'\sim z'$,
because $x'*y'\nsim u'*x'$, $x'*y'\nsim v'*y'$ and $T'$ is a
complete triad. Then $(u*x,x'),(v*y,y')$ and $(z,x'*y')$ are in
$\{\alpha,\beta\}^{\perp}$.

By a similar argument, if $x\nsim y$ and $u'\sim v'$ then
$(u,u'*x'),(v,v'*y')$ and $(u*v,z')$ are in
$\{\alpha,\beta\}^{\perp}$, where $\{u'*x',v'*y',z'\}$ is the
complete triad of $S'$ containing $u'*x'$ and $v'*y'$.

Now, assume that $x\nsim y$ and $u'\nsim v'$. If $u'=x'$ and
$v'=y'$, then $(a,a'),(b,b')$ and $(c,c')$ are in
$\{\alpha,\beta\}^{\perp}$, where $\{u',v'\}^{\perp}=\{a',b',c'\}$
in $S'$. We may assume that $v'\neq y'$. Then the complete triads
$\{x',y'\}^{\perp}$ and $\{u',v'\}^{\perp}$ of $S'$ intersect at
$w'=v'*y'$. This fact is independent of whether $u'=x'$ or not. Let
$\{x',y'\}^{\perp}=\{a',b',w'\}$ and
$\{u',v'\}^{\perp}=\{p',q',w'\}$ in $S'$. Since $a'\nsim w'$,
$b'\nsim w'$ and $\{p',q',w'\}$ is a complete triad of $S'$, each of
$a'$ and $b'$ is collinear with exactly one of $p'$ and $q'$.
Similarly, each of $p'$ and $q'$ is collinear with exactly one of
$a'$ and $b'$. So, we may assume that $a'\sim p'$ and $b'\sim q'$.
Then $(p,a'),(q,b')$ and $(w,w')$ are contained in
$\{\alpha,\beta\}^{\perp}$.
\end{proof}

\begin{lemma}\label{no-argument-coincide-distance-3}
Assume that $(A4)$ holds. Then $d(\alpha,\beta)= 3$.
\end{lemma}

\begin{proof}
We may assume that $u'\in y'^{\perp}$ and $v'\notin x'^{\perp}$.
Suppose that $d(\alpha,\beta)=2$ and $(z,w')\in
\{\alpha,\beta\}^{\perp}$. Then $z\notin\{x,y\}$,
$w'\notin\{u',v'\}$, $u',v'\in z'^{\perp}$ and
$w'\in\{x',y'\}^{\perp}$. Let $T'=\{x',y',z'\}$. Then $T'$ is either
a line or a complete triad of $S'$, because $x,y$ and $z$ are
pair-wise distinct and $u',w'\in T'^{\perp}$ with $u'\neq w'$. Since
$v'\in\{y',z'\}^{\perp}$, it follows that $v'\in T'^{\perp}$ and
$v'\in x'^{\perp}$, a contradiction to our assumption.

So $d(\alpha,\beta)\neq 2$. Now, choose $w'\in \{x',y'\}^{\perp}$
with $w'\neq u'$. Then $\alpha\sim (y,w')$ and $d((y,w'),\beta)=2$
by Lemma \ref{one-argument-coincide}. Hence $d(\alpha,\beta)= 3$.
\end{proof}

As a consequence of the above results, we have

\begin{corollary}\label{diameter-3}
The diameter of $\mathcal{S}$ is 3.
\end{corollary}

We next prove that the near-polygon property $(NP)$ is satisfied in
$\mathcal{S}$. Let $\mathtt{L}=\{\alpha,\beta,\gamma\}$ be a line
and $\theta$ be a point of $\mathcal{S}$. Let
$\alpha=(x,u'),\beta=(y,v'),\gamma=(z,w')$ and $\theta = (p,q')$.
Then $T=\{x,y,z\}$ is either a line or a complete triad of $S$ and
$T'^{\perp}=\{u',v',w'\}$.

\begin{proposition}\label{distance-2-to-two-points}
If $\theta$ has distance 2 from two points of $\mathtt{L}$, then it
is collinear with the third point of $\mathtt{L}$.
\end{proposition}

\begin{proof}
Let $d(\theta,\alpha)=d(\theta,\beta)=2$. We prove $\theta
\sim\gamma$ by showing that $p\neq z,q'\neq w'$, $q'\in z'^{\perp}$
and $w'\in p'^{\perp}$.

If $p=x$ and $q'\neq v'$ (respectively, $p\neq x$ and $q'=v'$), then
$d(\theta,\beta)=2$ (respectively, $d(\theta,\alpha)=2$) yields
$v'\notin x'^{\perp}$ (respectively, $u'\notin x'^{\perp}$), a
contradiction. So $p=x$ if and only if $q'=v'$. Similarly, $p=y$ if
and only if $q'=u'$. Thus, if $p\in\{x,y\}$, then $p\neq z,q'\neq
w'$, $q'\in z'^{\perp}$ and $w'\in p'^{\perp}$.

If $p\notin\{x,y\}$, then the above argument implies that $p\neq z$
and $q'\neq w'$. Also, $d(\theta,\alpha)=d(\theta,\beta)=2$ yields
$x',y'\notin q'^{\perp}$ and $u',v'\notin p'^{\perp}$. This implies
that $q'\in z'^{\perp}$ and $w'\in p'^{\perp}$.
\end{proof}

\begin{proposition}\label{distance-3-to-two-points}
If $\theta$ has distance 3 from two points of $\mathtt{L}$, then it
has distance 2 to the third point of $\mathtt{L}$.
\end{proposition}

\begin{proof}
Let $d(\theta,\alpha)=d(\theta,\beta)=3$. We prove
$d(\theta,\gamma)=2$. By Lemma \ref{one-argument-coincide}, we may
assume that $p\neq z$ and $q'\neq w'$. This together with
$d(\theta,\alpha)=d(\theta,\beta)=3$ imply that $p'\notin T'$,
$q'\notin T'^{\perp}$. We show that $q'\notin z'^{\perp}$ and
$w'\notin p'^{\perp}$. This would complete the proof.

Suppose that $q'\in z'^{\perp}$. Since $q'\notin T'^{\perp}$,
$q'\notin x'^{\perp}$ and $q'\notin y'^{\perp}$. Then,
$d(\theta,\alpha)=d(\theta,\beta)=3$ yields $u',v'\in p'^{\perp}$.
This implies that $p'\in \{u',v'\}^{\perp}=T'$, a contradiction. A
similar argument shows that if $w'\notin p'^{\perp}$, then $q'\in
T'^{\perp}$, a contradiction.
\end{proof}

\textbf{Proof of Theorem \ref{first-main-theorem}.} Propositions
\ref{distance-2-to-two-points} and \ref{distance-3-to-two-points}
together with Corollary \ref{diameter-3} imply that $\mathcal{S}$ is
a near hexagon. By Lemmas \ref{one-argument-coincide} and
\ref{no-argument-coincide-distance-2}, $\mathcal{S}$ is dense. Since
$|\mathcal{P}|=105$, Theorem \ref{classification-result} completes
the proof.\medskip

Thus, quads in $\mathcal{S}$ are $(2,1)$ or $(2,2)$-GQs. In fact,
it can be shown that equality holds in Lemmas
\ref{one-argument-coincide} and
\ref{no-argument-coincide-distance-2}.

\section{Proof of Theorem \ref{second-main-theorem}}
\label{Proof-of-second-main-theorem}

By the construction of lines of $\mathbb{S}$, no two points of $P$,
as well as of $P'$, are collinear in $\mathbb{S}$. Further, if $x\in
P$ and $u'\in P'$, then $x\sim u'$ if and only if $(x,u')\in
\mathcal{P}$, or equivalently, $u'\in x'^{\perp}$ in $S'$. Let
$\alpha$ and $\beta$ be two distinct non-collinear points of
$\mathbb{S}$ with $\alpha \in P\cup P'$. Then one of the following
possibilities occur:
\begin{enumerate}
\item[$(B1)$] $\alpha =x$ and $\beta =y$ for some $x,y\in P$ with
$x\neq y$;

\item[$(B2)$] $\alpha =u'$ and $\beta =v'$ for some $u',v'\in P'$
with $u'\neq v'$;

\item[$(B3)$] $\alpha =x\in P$ and $\beta =u'\in P'$ with
$u'\notin x'^{\perp}$;

\item[$(B4)$] $\alpha =x\in P$ and $\beta =(y,v')\in \mathcal{P}$
with $x\neq y$ and $v'\in x'^{\perp}$ in $S'$;

\item[$(B5)$] $\alpha =u'\in P'$ and $\beta =(y,v')\in
\mathcal{P}$ with $u'\neq v'$ and $y\in u^{\perp}$ in $S$;

\item[$(B6)$] $\alpha =x\in P$ and $\beta =(y,v')\in \mathcal{P}$
with $x\neq y$ and $v'\notin x'^{\perp}$ in $S'$;

\item[$(B7)$] $\alpha =u'\in P'$ and $\beta =(y,v')\in
\mathcal{P}$ with $u'\neq v'$ and $y\notin u^{\perp}$ in $S$.
\end{enumerate}

\begin{lemma}\label{distanct-x,y-in-P-u',v'-in-P'}
Assume that $(B1)$ or $(B2)$ holds. Then
$|\{\alpha,\beta\}^{\perp}|\geq 3$ in $\mathbb{S}$.
\end{lemma}

\begin{proof}
If $(B1)$ holds, then $w'\in\{x,y\}^{\perp}$ in $\mathbb{S}$ for
each $w'\in \{x',y'\}^{\perp}$ in $S'$. So
$|\{\alpha,\beta\}^{\perp}|\geq 3$. Similarly, if $(B2)$ holds
then $|\{\alpha,\beta\}^{\perp}|\geq 3$.
\end{proof}

\begin{lemma}\label{distance-between-x,u'}
Assume that $(B3)$ holds. Then $d(\alpha,\beta)=3$.
\end{lemma}

\begin{proof}
Clearly $d(\alpha,\beta)\geq 3$ since $u'\notin x'^{\perp}$. Let
$v'\in \{u',x'\}^{\perp}$ in $S'$. Then $x,v',v,u'$ is a path of
length 3 in $\Gamma(\mathbb{P})$. So $d(\alpha,\beta)=3$.
\end{proof}

\begin{lemma}\label{distance-2-x,(y,v')}
Assume that $(B4)$ or $(B5)$ holds. Then
$|\{\alpha,\beta\}^{\perp}|\geq 3$.
\end{lemma}

\begin{proof}
Assume that $(B4)$ holds. If $x\sim y$ in $S$, then $v'\in x'y'$ in
$S'$. We may assume that $v'\neq x'$. Then $v',(x,x')$ and
$(x,v'*x')$ are in $\{\alpha,\beta\}^{\perp}$. If $x\nsim y$ in $S$,
let $\{x',y'\}^{\perp}=\{u',v',w'\}$ in $S'$. Then $v',(x,u')$ and
$(x,w')$ are in $\{\alpha,\beta\}^{\perp}$. A similar argument holds
if $(B5)$ holds.
\end{proof}

\begin{lemma}\label{distance-3-x,(y,v')}
Assume that $(B6)$ or $(B7)$ holds. Then $d(\alpha,\beta)=3$.
\end{lemma}

\begin{proof}
Assume that $(B6)$ holds. Suppose that
$\theta\in\{\alpha,\beta\}^{\perp}$. Then $\theta\neq v'$, since
$v'\notin x'^{\perp}$ in $S'$. So $\theta=(x,w')$ for some $w'\in
x'^{\perp}$. Then $\theta\sim\beta$ implies that $v'\in x'^{\perp}$
in $S'$, a contradiction. So $d(\alpha,\beta)\neq 2$. Now, $y\sim
\beta$ and $d(\alpha,y)=2$ (Lemma
\ref{distanct-x,y-in-P-u',v'-in-P'}). So $d(\alpha,\beta)=3$. A
similar argument can be applied if $(B7)$ holds.
\end{proof}

As a consequence of the above results of this section together with
Corollary \ref{diameter-3}, we have

\begin{corollary}\label{diameter-3-second}
The diameter of $\mathbb{S}$ is 3.
\end{corollary}

Next we prove that the property $(NP)$ is satisfied in
$\mathbb{S}$.

\begin{proposition}\label{distance-2-to-two-points-second}
Let $\mathtt{L}$ be a line of $\mathbb{S}$ of type
$(\mathbb{L}_{1})$ and $\alpha$ be a point of $\mathbb{S}$ not
contained in $\mathtt{L}$. Then $\alpha$ is nearest to exactly one
point of $\mathtt{L}$.
\end{proposition}

\begin{proof}
Let $\mathtt{L}=\{x,\beta,u'\}$ where $\beta =(x,u')\in
\mathcal{P}$. Let $\alpha =v'\in P'$. Then $v'\neq u'$ and
$d(\alpha, u')=2$ (Lemma \ref{distanct-x,y-in-P-u',v'-in-P'}). Now
$d(\alpha,\beta)=2$ or 3 according as $x\in v^{\perp}$ in $S$ or
not. In the first case, $\alpha\sim x$, and in the later case,
$d(\alpha,x)=3$ (Lemma \ref{distance-2-x,(y,v')}). A similar
argument holds if $\alpha \in P$.

Let $\alpha =(y,v')\in \mathcal{P}$. If $x=y$, then $u'\neq v'$ and
$x\in\{u,v\}^{\perp}$ in $S$. So $\alpha\sim x$ and $d(\alpha,
\beta)=d(\alpha,u')=2$ (Lemmas \ref{one-argument-coincide} and
\ref{distance-2-x,(y,v')}). Similarly, if $u'=v'$ then $\alpha\sim
u'$ and $d(\alpha, \beta)=d(\alpha,x)=2$ . Assume that $x\neq y$ and
$u'\neq v'$. If $\alpha \sim \beta$, then $u'\in y'^{\perp}$ and
$v'\in x'^{\perp}$ in $S'$. So $d(\alpha,x)=d(\alpha,u')=2$ (Lemma
\ref{distance-2-x,(y,v')}). If $d(\alpha,\beta)=2$, then $u'\notin
y'^{\perp}$ and $v'\notin x'^{\perp}$ in $S'$. By Lemma
\ref{distance-3-x,(y,v')}, $d(\alpha,x)=d(\alpha,u')=3$. If
$d(\alpha,\beta)=3$, then either $u'\in y'^{\perp}$ and $v'\notin
x'^{\perp}$, or $u'\notin y'^{\perp}$ and $v'\in x'^{\perp}$ in
$S'$. Then $d(\alpha,x)=3$ and $d(\alpha,u')=2$ in the first case,
and $d(\alpha,x)=2$ and $d(\alpha,u')=3$ in the later.
\end{proof}

Now, let $\mathtt{L}=\{\beta,\theta,\gamma\}\in \mathcal{L}$ be a
line of $\mathbb{S}$ and $\alpha \in P\cup P'$. We take
$\beta=(x,u'),\theta=(y,v')$ and $\gamma =(z,w')$.

\begin{proposition}\label{distance-3-to-two-points-second}
If $\alpha$ has distance 2 from two points of $\mathtt{L}$, then it
is collinear with the third point of $\mathtt{L}$.
\end{proposition}

\begin{proof}
Let $\alpha=q' \in P'$ and $d(\alpha,\beta)=d(\alpha,\theta)=2$.
Then $q'\notin\{u',v'\}$ and $x,y\in q^{\perp}$ in $S$. Thus,
$q'\in\{x',y'\}^{\perp}=\{u',v',w'\}$ in $S'$. So $q'=w'$ and
$\alpha\sim\gamma$. A similar argument holds if $\alpha\in P$.
\end{proof}

\begin{proposition}\label{distance-3-to-two-points-second}
If $\alpha$ has distance 3 from two points of $\mathtt{L}$, then it
has distance 2 to the third point of $\mathtt{L}$.
\end{proposition}

\begin{proof}
Let $\alpha=q' \in P'$ and $d(\alpha,\beta)=d(\alpha,\theta)=3$.
Then $q'\notin\{u',v'\}$ and $x,y\notin q^{\perp}$ in $S$. So
$q'\neq w'$ and $q'\in z'^{\perp}$ in $S'$. The later follows
because $\{x,y,z\}$ is a line or a complete triad of $S$. Thus,
$d(\alpha,\gamma)=2$. A similar argument holds if $\alpha\in P$.
\end{proof}

\textbf{Proof of Theorem \ref{second-main-theorem}.} By the results
of this section together with Theorem \ref{first-main-theorem},
$\mathbb{S}$ is a slim dense near hexagon. Since
$|\mathcal{P}|=135$, Theorem \ref{classification-result} completes
the proof.

\end{document}